\documentclass[12 pt]{article}
\usepackage{times}
\usepackage[english]{babel}
\usepackage[T1]{fontenc}
\usepackage[latin1]{inputenc}
\usepackage{amsmath}
\usepackage{amsthm}
\usepackage{amsfonts}
\usepackage{eucal}

\newtheorem{lem}{Lemma}
\newtheorem{pro}{Proposition}
\newtheorem{rem}{Remark}
\newtheorem{teo}{Theorem}

\newcommand\F{{\mathbb F}}
\renewcommand\O{{\cal O}}
\newcommand\Z{{\mathbb Z}}
\newcommand\Q{{\mathbb Q}}
\newcommand\N{{\mathbb N}}
\newcommand\DIM{{\smallskip\noindent{\bf Proof.}\quad}}
\newcommand\CVD{\begin{flushright}$\square$\end{flushright}
\vskip 0.2cm\goodbreak}

\begin{document}

\title{On the linear independence of $p$-adic $L$-functions modulo $p$}
\date{ }
\maketitle
\begin{center}
\scshape{Bruno Anglès, Gabriele Ranieri}
\end{center}
\vskip 1cm
\begin{center}  
Laboratoire de math\'ematiques Nicolas Oresme,\\ 
CNRS UMR 6139, Universit\'e de Caen BP 5186,\\ 
14032 Caen cedex, France\\
E-mail: bruno.angles@math.unicaen.fr\\
E-mail: ranieri@math.unicaen.fr
\end{center}
\vskip 1cm

\begin{center}
{\scshape Abstract}
\end{center}
\noindent
Let $p \geq 3$ be a prime. 
Let $n \in \N$ such that $n \geq 1$, let $\chi_1, \ldots , \chi_n$ be characters of conductor $d$ not divided by $p$ and let $\omega$ be the Teichm\"uller character.  
For all $i$ between $1$ and $n$, for all $j$ such that $0 \leq j \leq ( p-3 )/2$, set 
$$
\theta_{i, j} =
\begin{cases}
\chi_i \omega^{2 j + 1}&\;\hbox{if} \ \chi_i\;{\rm is} \ {\rm odd}; \\
\chi_i \omega^{2 j}&\;\hbox{if} \ \chi_i\;{\rm is} \ {\rm even}.
\end{cases}
$$ 
Let $K$ be the least extension of $\Q_p$ which contains all the values of the characters $\chi_i$ and let $\pi$ be a prime of the valuation ring $\O_K$ of $K$. 
For all $i, j$ let $f ( T, \theta_{i, j} )$ be the Iwasawa series associated to $\theta_{i, j}$ and $\overline{f ( T, \theta_{i, j} )}$ its reduction modulo $( \pi )$.
Finally let $\overline{\F_p}$ be an algebraic closure of $\F_p$.
Our main result is that if the characters $\chi_i$ are all distinct modulo $( \pi )$, then $1$ and the series $\overline{f ( T, \theta_{i, j} )}$, for $1 \leq i \leq n$ and $0 \leq j \leq ( p- 3 )/2$, are linearly independent over a certain field $\Omega$ which contains $\overline{\F_p} ( T )$.      

\section{Introduction}\label{sec:parteprima}

Let $ p $ be an odd prime.
Let $n \in \N$ such that $n \geq 1$, let $\chi_1, \ldots , \chi_n$ be characters of conductor $ d $ not divided by $p$ and let $\omega$ be the Teichm\"uller character.  
For all $ i $ between $1$ and $ n $ and $j$ such that $0 \leq j \leq ( p-3 )/2$, set
$$
\theta_{i, j} =
\begin{cases}
\chi_i \omega^{2 j + 1}&\;\hbox{if} \ \chi_i\;{\rm is} \ {\rm odd}; \\
\chi_i \omega^{2 j}&\;\hbox{if} \ \chi_i\;{\rm is} \ {\rm even}.
\end{cases}
$$  
Observe that, by definition, $ \theta_{i, j} $ is an even character for all $i, j$.

Set $\kappa_0 = 1 + dp$ and $K = \Q_p ( \chi_1, \ldots , \chi_n )$ (i.e. the least extension of $\Q_p$ generated by all the images of $\chi_i$ for all $i$).
Let $\pi$ be a prime of the valuation ring $\O_K$ of $ K $ and let $\F_q$ be $\O_K/ \pi \O_K$. 
For all $ i $ between $1$ and $ n $ and $j$ between $0$ and $( p-3 )/ 2$, set $f ( T, \theta_{i, j} )$ the Iwasawa power series attached to the $p$-adic $L$-function $L_p ( s, \theta_{i, j} )$ (see~\cite[Theorem 7.10]{4}) and $\overline{f ( T, \theta_{i, j} )}$ its reduction modulo $( \pi )$.

Let $\overline{\F_p}$ be an algebraic closure of $\F_p$. 
For $F ( T ) \in \overline{\F_p}[[T]]$, we say that $ F ( T ) $ is a pseudo-polynomial if and only if there exist $r \in \N$, $a_1, \ldots , a_r \in \Z_p$ and $c_1, \ldots , c_r \in \overline{\F_p}$ such that
\[
F ( T ) = \sum_{i = 1}^r c_i ( T + 1 )^{a_i}.
\]
Then the set of the pseudo-polynomials is a ring which we denote by $A$. 
Moreover we denote by $ \Omega $ the quotient field of $ A $. 
The elements of $ \Omega $ are called pseudo-rational functions. 
Anglès (see~\cite[Theorem 4.5]{1}) shows that for all non-trivial even character of the first kind $ \theta $, $\overline{f ( T, \theta )}$ is not a pseudo-rational function. 
We shall prove the following generalisation of this result:

\begin{teo}\label{teo:teo1}
Suppose that the characters $ \chi_i $ are all distinct modulo $( \pi )$ (i.e. for all integer $i \neq j$ there exist $a \in ( \Z / d \Z )^\ast$ such that $\chi_i ( a ) \not \equiv \chi_j ( a ) \mod ( \pi )$). 
Then the elements of the set
\[
\{ 1, \ \overline{f ( T, \theta_{i, j} )}, \ 1 \leq i \leq n, \ 0 \leq j \leq ( p-3 )/ 2 \} 
\]
are linearly independent over $\Omega$.
\end{teo}

Observe that in the statement of Theorem~\ref{teo:teo1} it is necessary to suppose that the characters $\chi_i$ are all distinct modulo $( \pi )$.
Indeed suppose that there exist $i \neq k$ between $1$ and $n$ such that $\chi_i$ and $ \chi_k $ are congruent modulo $( \pi )$.
Since $p$ is odd this implies that $\chi_i$ and $\chi_k$ have the same parity.
Then for all $0 \leq j \leq ( p-3 )/2$ we have that $\theta_{i, j}$ is congruent to $\theta_{k, j}$ modulo $( \pi )$, which implies $\overline{f ( T, \theta_{i, j} )} = \overline{f ( T, \theta_{k, j} )}$.
Thus in this case the series $\overline{f ( T, \theta_{i, j} )}$ are dependent.

Observe also that if the characters $\chi_i$ are distinct modulo $( \pi )$, then for all $i, i^\prime$ between $1$ and $n$, $j, j^\prime$ between $0$ and $( p-3 )/2$, $\theta_{i, j}$ is congruent to $\theta_{i^\prime, j^\prime}$ modulo $( \pi )$ if and only if $i = i^\prime, j = j^\prime$.
It is clear that $i = i^\prime$ implies $j = j^\prime$ and $j = j^\prime$ implies $i = i^\prime$.  
Then suppose that $i \neq i^\prime$, $j \neq j^\prime$ and that $\theta_{i, j}$ is congruent to $\theta_{i^\prime, j^\prime}$ modulo $( \pi )$. 
Moreover suppose that $ \chi_i $ and $\chi_{i^\prime}$ are even (the other case is identical). 
Hence there exists an integer $ a $ such that $\omega^{2j} ( a ) \not \equiv \omega^{2j^\prime} ( a ) \mod ( \pi )$.
Since $p$ does not divide $ d $ there exists an integer $c$ such that $1 + cd \equiv a \mod ( p )$.
Then $\theta_{i, j} ( 1 + cd ) \equiv \omega^{2j} ( a ) \mod ( \pi )$ and $\theta_{i^\prime, j^\prime} ( 1+cd ) \equiv \omega^{2j^\prime} ( a ) \mod ( \pi )$.
Since $\theta_{i, j}$ and $\theta_{i^\prime, j^\prime}$ are equivalent modulo $( \pi )$, we get $\omega^{2j} ( a ) \equiv \omega^{2j^\prime} ( a ) \mod ( \pi )$ which is a contradiction.      

As in the proof of~\cite[Theorem 4.5]{1}, the main ingredient in the proof of Theorem~\ref{teo:teo1} is a remarkable result due to Sinnott.
Before the statement of that result we must define the following equivalence relation: let $a, b \in \Z_p^\ast$.
We say that $a$ is equivalent to $b \mod ( \Q^\ast )$ ($a \equiv b \mod ( \Q^\ast )$) if and only if there exists $c \in \Q^\ast$ such that $a b^{-1} = c$.

\begin{pro}\label{pro:pro2}{\rm (\cite[Proposition 1]{3})}
Let $ F $ be a finite field of characteristic $p$ and let $r_1 ( T ), \ldots , r_s ( T ) \in F ( T ) \cap F[[T]]$.
Let $c_1, \ldots , c_s \in \Z_p- \{ 0 \}$ and suppose that
\[
\sum_{i = 1}^s r_i ( ( T + 1 )^{c_i} - 1 ) = 0.
\]         
Then for all $a \in \Z_p$,
\[
\sum_{c_i \equiv a \mod ( \Q^\ast )} r_i ( ( T + 1 )^{c_i} - 1 ) \in F.
\]      
\end{pro}

Let's describe briefly the strategy of the proof of Theorem~\ref{teo:teo1}.
First (section~\ref{sec:parteseconda}) we recall some properties of the $p$-adic Leopoldt transform (most of them already proved in~\cite{1}) which we will often use.         

Then (section~\ref{sec:parteterza}) we consider the case $d \geq 2$.
Some results about the $p$-adic Leopoldt transform of~\cite{1} and Proposition~\ref{pro:pro2} will allow us to reduce the proof of Theorem~\ref{teo:teo1} to the computation of the rank of a certain matrix whose entries depend on the values of the characters $ \chi_i $ (Lemma~\ref{lem:lem11}). 
After such computation the proof of this case of the theorem will follow by some simple remarks of linear algebra.

In section~\ref{sec:partequarta} we study the case $d = 1$.
In that case we have to consider a <<perturbation>> of the functions $f ( T, \theta_{i, j} )$ to be able to apply Proposition~\ref{pro:pro2}.
Then the proof is not very different from the proof of the previous case (actually it is simpler since it does not request a result similar to Lemma~\ref{lem:lem11}) and some remarks of linear algebra will imply the assumption.                     

Finally we give a possible link between Theorem~\ref{teo:teo1} and Ferrero-Washington's heuristic (see~\cite{2}).  
Let $i$ be an integer between $1$ and $( p-3 )/2$.
Write   
\[
f ( T, \omega^{2i} ) = \sum_{k = 0}^{+ \infty} a_k ( \omega^{2i} ) T^k.
\]
The $\lambda$-invariant of $f ( T, \omega^{2i} )$, denoted by $\lambda ( \omega^{2i} )$, is the least $ k $ such that $a_k ( \theta ) \not \equiv 0 \mod ( p )$.
We set 
\[
\lambda^- = \sum_{i = 1}^{( p-3 ) /2} \lambda ( \omega^{2i} ).
\]
Ferrero and Washington make the following hypothesis to define a heuristic to make previsions about possible bounds for $\lambda^-$.   
\vspace{5 pt}

{\it
Ferrero-Washington's hypothesis: Every coefficient of $ f ( T, \omega^{2i} ) $ is random$\mod ( p )$ and independent from the other coefficients.
}
\vspace{5 pt}

Theorem~\ref{teo:teo1} implies that  
\[
1, \overline{f ( T, \omega^0 )}, \overline{f ( T, \omega^2 )}, \ldots ,\overline{f ( T, \omega^{p-3} )}    
\]
are linearly independent over $\Omega$. 
Thus our result seems to confirm Washington's hypothesis. 
 
\section{Preliminaries}\label{sec:parteseconda} 
 
In this section we shall list some properties of the $p$-adic Leopoldt transform which will be very important in the proof of Theorem~\ref{teo:teo1}.
Let $L$ be a finite extension of $\Q_p$, $\O_L$ its valuation ring and $\F_{q^\prime}$ its residue field.
Let $\kappa$ a topological generator of $1 + p \Z_p$ and, for all $a \in \Z_p^\ast$, set $\omega ( a )$ the unique $( p-1 )$th root of unity in $ \Z_p $ congruent to $a \mod ( p )$.     
Following~\cite{1} for all $\delta \in \Z / ( p -1 ) \Z$ we define $p$-adic Leopoldt transform $ \Gamma_\delta $ the unique continuous $\mathcal{O}_L$-linear endomorphism of $\mathcal{O}_L[[T]]$ such that for all $a \in \Z_p$,
$$
\Gamma_\delta ( ( T + 1 )^a ) = 
\begin{cases}
\omega^\delta ( a ) ( T + 1 )^{{\log_p ( a ) \over \log_p ( \kappa )}}&\;\hbox{if} \ a \in \Z_p^\ast; \\
0&\;\hbox{otherwise} 
\end{cases}
$$    
(see~\cite[Sections 2., 3.]{1} for the proof of the fact that $\Gamma_\delta$ is well-defined and unique).
In an obvious way we can define a similar $\F_{q^\prime}$-linear continuous endomorphism of $\F_{q^\prime}[[T]]$ that we denote by $\overline{\Gamma}_\delta$.
Observe that if $a \in \Z_p^\ast$ we have $a \equiv \omega ( a ) \mod ( p )$.
Thus, for all $a \in \Z_p^\ast$, we have
\[
\overline{\Gamma}_\delta ( ( T + 1 )^a ) = a^\delta ( T + 1 )^{{\log_p ( a ) \over \log_p ( \kappa )}}.
\]

In the proof of Theorem~\ref{teo:teo1} we use other $\O_L$-linear endomorphisms of $\O_L[[T]]$ already introduced by Anglès in~\cite{1}. 
Let us recall their definition. 
Let $\mu_{p-1} \subseteq \Z_p^\ast$ be the group of $( p-1 )$th roots of unity. 
For all $\delta \in \Z/ ( p-1 ) \Z$ and $F ( T ) \in \O_L[[T]]$, we set
\[
\gamma_\delta ( F ( T ) ) = {1 \over p-1} \sum_{\eta \in \mu_{p-1}} \eta^\delta F ( ( T + 1 )^\eta - 1 ).
\] 
Observe that $\gamma_\delta \gamma_{\delta^\prime} = 0$ for $\delta \neq \delta^\prime$, $\gamma_\delta^2 = \gamma_\delta$ and $\sum_{\delta \in \Z/ ( p-1 ) \Z} \gamma_\delta = Id_{\O_L[[T]]}$. 

For $F ( T ) \in \O_L[[T]]$ set
\begin{align*}
D ( F ( T ) ) & = ( T + 1 ) {d \over dT} F ( T ) \\
U ( F ( T ) ) & = F ( T ) - {1 \over p} \sum_{\zeta \in \mu_p} F ( \zeta ( T + 1 ) - 1 ) \in \O_L[[T]].
\end{align*}
In an obvious way we can define the $\F_{q^\prime}$-linear endomorphism of $\F_{q^\prime}[[T]]$ $\overline{\gamma}_\delta$, $\overline{D}$ and $\overline{U}$.    
Observe that:
\begin{itemize}
\item $U^2 = U$;
\item $DU = UD$;
\item $\gamma_\delta U = U \gamma_\delta$ for all $\delta \in \Z/ ( p-1 ) \Z$;
\item $D \gamma_\delta = \gamma_{\delta + 1} D$ for all $\delta \in \Z/ ( p-1 ) \Z$;
\item $\overline{U} = \overline{D}^{p-1}$. 
\end{itemize} 

In the following lemma we shall list some properties of $ \Gamma_\delta $ whose we need to prove Theorem~\ref{teo:teo1}.  

\begin{lem}\label{lem:lem3}
Let $\delta \in \Z/ ( p- 1 ) \Z$ and $F ( T ) \in \O_L[[T]]$.
Then
\begin{enumerate}
\item $\Gamma_\delta ( F ( T ) ) = \Gamma_\delta \gamma_{-\delta} ( F ( T ) ) = \Gamma_\delta \gamma_{-\delta} U F ( T )$.\\
\item Suppose that $\Gamma_\delta ( F ( T ) )$ is a pseudo-polynomial.
Then $\gamma_{-\delta} U ( F ( T ) )$ is a pseudo-polynomial.\\ 
\item $\overline{\Gamma}_{\delta + 1} ( \overline{F ( T )} ) = \overline{\Gamma}_\delta \overline{D} ( \overline{F ( T )} )$.
\end{enumerate} 
\end{lem}

\DIM 1. See~\cite[Proposition 3.2(2)]{1}.

2. The assumption immediately follows from~\cite[Proposition 3.1]{1}.

3. Since $\overline{\Gamma}_\delta$ is a $\F_{q^\prime}$-linear continuous endomorphism of $\F_{q^\prime}[[T]]$, it suffices to prove that the assumption is true for $F ( T ) = ( T + 1 )^a$ with $a \in \Z_p$.
If $ p $ divides $a$ we have
\[
0 = \overline{\Gamma}_{\delta + 1} ( ( T + 1 )^a ) = \overline{\Gamma}_\delta ( \overline{D} ( ( T + 1 )^a ) )    
\]
and the assumption is trivial.

Suppose that $p$ does not divide $a$.
We have
\begin{align*}
\overline{\Gamma}_\delta (\overline{D} ( ( T + 1 )^a ) ) & = \overline{\Gamma}_\delta ( a ( T + 1 )^a ) )\\
& = a \overline{\Gamma}_\delta ( ( T + 1 )^a ) \\
& = a^{\delta + 1} ( T + 1 )^{{\log_p ( a ) \over \log_p ( \kappa )}}\\
& = \overline{\Gamma}_{\delta + 1} ( ( T + 1 )^a ).
\end{align*}
\CVD       

Let $\theta$ be a Dirichlet character of the first kind such that 
\[
\theta = \chi \omega^{\delta + 1},
\]
with $\delta \in \Z/ ( p-1 ) \Z$ and $\chi$ a character of conductor $d$ not divided by $p$.
Observe that $\kappa_0 = 1 + dp$ is a topological generator of $1 + p \Z_p$ and from now on set $\kappa = \kappa_0$. 
Suppose that $\Q_p ( \chi ) \subseteq L$.
Set
\[
F_\chi ( T ) = \sum_{a = 1}^d {\chi ( a ) ( T + 1 )^a \over 1 - ( T + 1 )^d}
\]
and $\overline{F_\chi} ( T ) \in \F_{q^\prime} ( T )$ its reduction modulo the maximal ideal of $\O_L$.
In the following lemma we list some properties of $F_\chi ( T )$ and we recall the relation between $F_\chi ( T )$ and $f ( T, \theta )$.

\begin{lem}\label{lem:lem4}
We have:
\begin{enumerate}
\item If $d \geq 2$, then $F_\chi ( T ) \in \O_L[[T]]$.\\
\item If $d = 1$, then $\gamma_\alpha ( F_\chi ( T ) ) \in \O_L[[T]]$ for all $\alpha \in \Z/ ( p-1 ) \Z$ and $\alpha \neq 1$.\\
\item If $d \geq 2$, then $F_\chi ( ( T + 1 )^{-1} - 1 ) = \varepsilon F_\chi ( T )$ where $\varepsilon = 1$ if $\chi$ is odd and $\varepsilon = -1$ if $\chi$ is even.\\
\item If $d = 1$, then $F_\chi ( ( T + 1 )^{-1} - 1 ) = - F_\chi ( T ) - 1$.\\  
\item If $d$ divides the positive integer $g$ we have
\[
F_\chi ( T ) = {\sum_{a = 1}^g \chi ( a ) ( T + 1 )^a \over 1 - ( T + 1 )^g}.
\]       
\item $\Gamma_\delta \gamma_{-\delta} U ( F_\chi ( T ) ) = f ( ( T + 1 )^{-1} - 1, \theta )$.
\end{enumerate} 
\end{lem}

\DIM For 1., 2., 3., 4. see~\cite[Lemma 4.1]{1}.

5. We have: 
\begin{align*}
{\sum_{a = 1}^g \chi ( a ) ( T + 1 )^a \over 1 - ( T + 1 )^g} & = \sum_{a = 1}^d \chi ( a ) \sum_{b \equiv a \mod ( d )} {\chi ( b ) ( T + 1 )^b \over 1 - ( T + 1 )^g} \\ 
& = \sum_{a = 1}^d \chi ( a ) { ( T + 1 )^a + ( T + 1 )^{a + d} + \ldots + ( T + 1 )^{a + g - d} \over 1 - ( T + 1 )^g} \\
& = {\sum_{a = 1}^d \chi ( a ) ( T + 1 )^a ( 1 - ( T + 1 )^g ) \over ( 1 - ( T + 1 )^g ) ( 1 - ( T + 1 )^d )} \\
& = F_\chi ( T ).
\end{align*}

6. Remember that $U \gamma_{-\delta} = \gamma_{-\delta} U$. 
Then apply~\cite[Lemma 4.4]{1}.  
\CVD

Consider the ring $\O_L[[t]]$. 
It acts over $\O_L[[T]]$ via: 
\[
( t + 1 ) F ( T ) = F ( ( T + 1 )^{\kappa_0} - 1 ) \in \O_L[[T]],
\]
for all $F ( T ) \in \O_L[[T]]$.

\begin{lem}\label{lem:lem5} 
Let $H ( T ) \in \O_L[[T]]$.
Then for all $F ( T ) \in \O_L[[T]]$ we have
\[
H ( T ) \Gamma_\delta ( F ( T ) ) = \Gamma_\delta ( H ( t ) F ( T ) ).
\]
\end{lem}

\DIM Since $\Gamma_\delta$ is $\O_L$-linear and continuous it suffices to prove the assumption in the case $P ( T ) = ( T + 1 )^a$, where $a \in \Z_p$.
By \cite[Proposition 3.2(3)]{1} for all $b \in \Z_p^\ast$ we have
\[
\Gamma_\delta ( F ( T + 1 )^b - 1 ) = \omega^\delta ( b ) ( T + 1 )^{{\log_p ( b ) \over \log_p ( \kappa )}} \Gamma_\delta ( F ( T ) ).
\] 
Then we get
\begin{align*}
\Gamma_\delta ( ( t + 1 )^a F ( T ) ) & = \Gamma_\delta ( F ( ( T + 1 )^{\kappa_0^a} - 1 )\\
& = \omega^\delta ( a ) ( T + 1 )^{{\log_p ( \kappa_0^a ) \over \log_p ( \kappa_0 )}} \Gamma_\delta ( F ( T ) )\\
& = ( T + 1 )^a \Gamma_\delta ( F ( T ) ).
\end{align*}
\CVD     
         
\section{The case $d \geq 2$}\label{sec:parteterza} 

The aim of this section is to prove Theorem~\ref{teo:teo1} in the case where the conductor $d$ of the characters $\chi_i$ is $\geq 2$.

{\bf Proof of Theorem~\ref{teo:teo1} in the case $d \geq 2$.} 
Let $\chi_1, \chi_2, \ldots , \chi_n$ characters of conductor $d \geq 2$ distinct modulo $( \pi )$.
Without loss of generality we can suppose $\chi_1, \ldots , \chi_r$ odd and $\chi_{r + 1}, \ldots , \chi_n$ even for a certain integer $r \leq n$.
For all $i$ between $1$ and $n$ and $j$ between $0$ and $( p-3 )/2$ set   
\begin{equation}\label{eqn:relcollina}
\theta_{i, j} =
\begin{cases}
\chi_i \omega^{2 j + 1}&\;\hbox{if} \ 1 \leq i \leq r; \\
\chi_i \omega^{2 j}&\;\hbox{otherwise.} 
\end{cases}
\end{equation}
Suppose that for all $1 \leq i \leq n$, $0 \leq j \leq ( p-3 )/2$, there exist $g_{i, j} ( T ) \in \Omega$ such that 
\[
\sum_{i = 1}^n \sum_{j = 0}^{( p- 3 )/2} g_{i, j}( T ) \overline{f( T, \theta_{i, j} )} \in \Omega
\] 
and $g_{i, j} ( T ) \neq 0$ for some $i, j$.
Set $h_{i, j} ( T ) = g_{i, j} ( 1/ ( T + 1 ) - 1 )$ for all $i, j$. 
Then we have
\[
\sum_{i = 1}^n \sum_{j = 0}^{( p-3 )/ 2} h_{i, j}( T ) \overline{f( 1/ ( T + 1 ) -1, \theta_{i, j} )} \in \Omega
\]
and $h_{i, j} ( T ) \neq 0$ for certain $i, j$. 
Observe that we can suppose that $h_{i, j} \in A$ for all $i, j$ and that 
\begin{equation}\label{eqn:rel41}
\sum_{i = 1}^n \sum_{j = 0}^{( p-3 )/ 2} h_{i, j}( T ) \overline{f( 1/ ( T + 1 ) -1, \theta_{i, j} )} \in A.
\end{equation}
By Lemma~\ref{lem:lem4}(6.) and~(\ref{eqn:relcollina}) for all $i$ between $1$ and $n$ and $j$ between $0$ and $( p-3 ) /2$, we have:
\begin{equation}\label{eqn:relcollina2}
f\bigg( {1 \over T + 1} - 1, \theta_{i, j} \bigg) = 
\begin{cases}
\Gamma_{2j} \gamma_{-2j} ( F_{\chi_i} ( T ) )&\;\hbox{if} \ 1 \leq i \leq r; \\
\Gamma_{2j - 1} \gamma_{-2j + 1}( F_{\chi_i} ( T ) )&\;\hbox{otherwise.} 
\end{cases}
\end{equation}
Moreover by Lemma~\ref{lem:lem3}(3.), for all $\delta \in \Z/ ( p- 1 )\Z$ we have 
\[
\overline{\Gamma}_{\delta + 1} \overline{\gamma}_{- \delta - 1} = \overline{\Gamma}_\delta \overline{\gamma}_\delta \overline{D}.  
\]
Hence we can rewrite relation~(\ref{eqn:rel41}) in the following way: 
\[
\sum_{i = 1}^r \sum_{ j = 0}^{ ( p-3 ) /2 } h_{i, j}( T ) \overline{\Gamma}_0 \overline{\gamma}_0 ( \overline{D}^{2j} \overline{F}_{\chi_i} ( T ) ) + \sum_{i = r+ 1}^n \sum_{j = 0}^{ ( p-3 ) /2 }  h_{i, j}( T ) \overline{\Gamma}_0 \overline{\gamma}_0 ( \overline{D}^{2j - 1} \overline{F}_{\chi_i} ( T ) ) \in A.
\]  
By Lemma~\ref{lem:lem3}(1.) this last relation implies that 
\begin{equation}\label{eqn:rel42}
\sum_{i = 1}^r \sum_{ j = 0}^{( p- 3 )/2} h_{i, j}( T ) \overline{\Gamma}_0 \overline{\gamma}_0 \overline{U} ( \overline{D}^{2j} \overline{F}_{\chi_i} ( T ) ) + \sum_{i = r+ 1}^n \sum_{j = 0}^{( p - 3 )/2} h_{i, j}( T ) \overline{\Gamma}_0 \overline{\gamma}_0 \overline{U} ( \overline{D}^{2j - 1} \overline{F}_{\chi_i} ( T ) ) \in A.  
\end{equation}
Applying Lemma~\ref{lem:lem5} and~\ref{lem:lem3}(2.) to~(\ref{eqn:rel42}) we get
\begin{equation}\label{eqn:rel43}
\sum_{i = 1}^r \sum_{ j = 0}^{( p-3 )/2} h_{i, j}( t ) \overline{\gamma}_0 \overline{U} ( \overline{D}^{2j} \overline{F}_{\chi_i} ( T ) ) + \sum_{i = r+ 1}^n \sum_{j = 0}^{( p-3 )/2} h_{i, j} ( t ) \overline{\gamma}_0 \overline{U} ( \overline{D}^{2j - 1} \overline{F}_{\chi_i} ( T ) ) \in A.
\end{equation}     

Set for all $1 \leq i \leq n$, $0 \leq j \leq ( p-3 )/2$, 	 
\begin{equation}\label{eqn:relcollina3}
F_{i, j} ( T ) = 
\begin{cases}
\overline{U} \overline{D}^{2j} ( \overline{F}_{\chi_i} ( T ) )&\;\hbox{if} \ 1 \leq i \leq r; \\
\overline{U} \overline{D}^{2j - 1} ( \overline{F}_{\chi_i} ( T ) )&\;\hbox{otherwise.}
\end{cases}        
\end{equation}
Then we can rewrite relation~(\ref{eqn:rel43}) in the following way:
\begin{equation}\label{eqn:rel44}
\sum_{i = 1}^n \sum_{ j = 0}^{( p - 3 )/2} h_{i, j}( t ) \overline{\gamma}_0 ( F_{i, j} ( T ) ) \in A.
\end{equation}

Now recall that by Lemma~\ref{lem:lem4}(3.), we have
$$
F_{\chi_i} ( ( T + 1 )^{-1} - 1 ) =
\begin{cases}
F_{\chi_i} ( T )&\;\hbox{if} \ 1 \leq i \leq r ; \\
- F_{\chi_i} ( T )&\;\hbox{if} \ r+1 \leq i \leq n. 
\end{cases}
$$
Moreover observe that for all $1 \leq i \leq n$, $0 \leq k \leq p-2$ we have
\[
( \overline{U} \overline{D}^k \overline{F}_{\chi_i} ) ( ( T + 1 )^{-1} -1 ) = ( -1 )^k \overline{U} \overline{D}^k ( \overline{F}_{\chi_i} ( ( T + 1 )^{-1} -1 ) ).        
\] 
Let $i$ be between $1$ and $r$ and $j$ be between $0$ and $( p-3 )/ 2$.
Then 
\begin{align*}
F_{i, j} ( ( T + 1 )^{-1} -1 ) & = ( \overline{U} \overline{D}^{2j} \overline{F}_{\chi_i} ) ( ( T + 1 )^{-1} - 1 ) ) \\ 
& = (-1)^{2j} \overline{U} \overline{D}^{2j} ( \overline{F}_{\chi_i} ( ( T + 1 )^{-1} - 1 ) ) \\
& = \overline{U} \overline{D}^{2j} ( \overline{F}_{\chi_i} ( T ) ) \\
& = F_{i, j} ( T ).
\end{align*} 
With exactly the same computation we can prove that
\[
F_{i, j} ( ( T + 1 )^{-1} -1 ) = F_{i, j} ( T )  
\]
for all $i, j$, also in the case where $r + 1 \leq i \leq n$.
Then
\begin{equation}\label{eqn:rel45}    
F_{i, j} ( ( T + 1 )^{-1} - 1 ) = F_{i, j} ( T ), \ \forall i, j.   
\end{equation}

We recall that $\F_q$ is, by definition, the residue field of the least extension of $\Q_p$ which contains all the images of the characters $\chi_i$. 
Consider the least field which contains $\F_q$ and all the coefficients of $h_{i, j} ( T )$ for all $i, j$.
Since $h_{i, j} ( T ) \in A \subseteq \overline{F}_p[[T]]$, such field is a finite extension of $\F_q$.
Call it $\F_{q_1}$ and write $h_{i, j} ( t ) = \sum_{b \in \Z_p} c_{i, j, b} ( t + 1 )^b$ with $c_{i, j, b} \in \F_{q_1}$.
Moreover observe that since $h_{i, j} ( T ) \neq 0$ for certain integer $i, j$, there exist $i, j, b$ such that $c_{i, j, b} \neq 0$.
Let 
\[
G_b ( T ) = \sum_{i = 1}^n \sum_{ j = 0}^{( p- 3 )/2} c_{i, j, b} F_{i, j} ( T ).   
\]
By Lemma~\ref{lem:lem4}(1.), $G_b ( T ) \in \F_{q_1}[[T]] \cap \F_{q_1} ( T )$.
Since 
\[
( t + 1 )^b ( G_b ( T ) ) = G_b ( ( T + 1 )^{\kappa_0^b} - 1 ), 
\]
by~(\ref{eqn:rel44}) we have
\begin{equation}\label{eqn:rel46}
\overline{\gamma}_0 ( \sum_{b \in \Z_p}  G_b ( ( T + 1 )^{\kappa_0^b} - 1 ) ) \in A.
\end{equation}

Choose a subset of $\mu_{p-1}$ whose elements represent all the classes of $\mu_{p - 1} / \{-1, 1\}$ and call it $S$.
From~(\ref{eqn:rel45}) it follows that 
\begin{equation}\label{eqn:relpure}
G_b ( T ) = G_b ( ( T + 1 )^{-1} - 1 )
\end{equation} 
for all $b$.
Thus by~(\ref{eqn:rel46}) we get
\[
\sum_{\eta \in S} \sum_{b \in \Z_p}  G_b ( ( T + 1 )^{\eta \kappa_0^b} - 1 ) \in A.
\]
Since $G_b ( T ) = 0$ for all but finitely many $b \in \Z_p$, there exists a positive integer $u$ such that
\[
\sum_{\eta \in S} \sum_{k = 1}^u  G_{b_k} ( ( T + 1 )^{\eta \kappa_0^{b_k}} - 1 ) \in A.
\]
Moreover we set $G_k ( T ) = G_{b_k} ( T )$.     
By Proposition~\ref{pro:pro2} there exist an integer $l \leq u$, $b_1, b_2, \ldots b_l \in \Z_p$, $b_i \neq b_j$ for $i \neq j$, $\eta_1, \eta_2, \ldots \eta_l \in \mu_{p - 1}$ with $\eta_i \kappa_0^{b_i} \sim_{\Q^\ast} \eta_j \kappa_0^{b_j}$ for all $i$, $j$ and $\eta_i \kappa_0^{b_i} \neq \eta_j \kappa_0^{b_j}$ for $i \neq j$ such that   
\begin{equation}\label{eqn:rel47}
\sum_{k = 1}^l  G_k ( ( T + 1 )^{\eta_k \kappa_0^{b_k}} - 1 ) \in A.   
\end{equation}

For all $1 \leq k \leq l$ write
\[
\eta_k \kappa_0^{b_k} = \eta_1 \kappa_0^{b_1} x_k,
\]  
where $x_j \in \Q^\ast \cap \Z_p^\ast$ and $x_i \neq x_k$ if $i \neq k$.
Recall that by~(\ref{eqn:relpure}), we have $G_k ( T ) = G_k ( ( T + 1 )^{-1} - 1 )$.
Hence we can suppose that $x_k > 0$ for all $k$. 
By~(\ref{eqn:rel47}) we get
\[
\sum_{k = i}^l  G_k ( ( T + 1 )^{x_k} - 1 ) \in A.
\] 
Therefore there exist some positive integers $N_1 , N_2 , \ldots N_l$ not divided by $p$ such that $1 \leq N_1 < N_2 < \ldots < N_l$ and
\[
\sum_{k = i}^l  G_k ( ( T + 1 )^{N_k} - 1 ) \in A.    
\]
If we rewrite $G_k ( T )$ as a combination of $F_{i, j} ( T )$, we get 
\[
\sum_{k = 1}^l \sum_{i = 1}^n \sum_{ j = 0}^{( p - 3 )/2} c_{i, j, k} F_{i, j} ( ( T + 1 )^{N_k} - 1 ) \in A.  
\]
By \cite[Lemma 3.5]{1} if $H ( T ) \in \F_{q_1} ( T )$, then $H ( T ) \in A$ if and only if there exists $m \in \N$ such that $( T + 1 )^m H ( T ) \in \F_{q_1} [T]$.
Since the denominator of  $F_{i, j} ( ( T + 1 )^{N_k} - 1 )$ is relatively prime with $( T + 1 )$ for all $i, j, k$, we get
\begin{equation}\label{eqn:rel48}
\sum_{k = 1}^l \sum_{i = 1}^n \sum_{ j = 0}^{( p - 3 )/2} c_{i, j, k} F_{i, j} ( ( T + 1 )^{N_k} - 1 ) \in \F_{q_1} [T].  
\end{equation}

Observe that, by Lemma~\ref{lem:lem4}(5.), we have
\[
F_{\chi_i} ( T ) = \sum_{a = 1}^{d p} {\chi_i ( a ) ( T + 1 )^a \over 1 - ( T + 1 )^{d p }}.
\]           
Then for all $1 \leq i \leq n$ and $0 \leq j \leq ( p-3 )/2$, we have
\begin{equation}\label{eqn:relsola1}
F_{i, j} ( T ) = 
\begin{cases}
\sum_{a = 1, p \not \mid a}^{d p - 1} {\chi_i ( a ) a^{2j} ( T + 1 )^a \over 1 - ( T + 1 )^{d p }}&\;\hbox{if} \ 1 \leq i \leq r; \\ 
\sum_{a = 1, p \not \mid a}^{d p - 1} {\chi_i ( a ) a^{2j - 1} ( T + 1 )^a \over 1 - ( T + 1 )^{d p }} &\;\hbox{otherwise.}
\end{cases}      
\end{equation}

Then replacing in~(\ref{eqn:rel48}) using~(\ref{eqn:relsola1}), we get
\begin{equation}\label{eqn:rel49}
\sum_{k = 1}^l \sum_{i = 0}^r \sum_{j = 0}^{( p-3 )/2} c_{i, j, k} \sum_{a = 1 , p \not \mid a}^{dp} \bigg( {a^{2j} \chi_i ( a ) ( T + 1 )^{a N_k} \over 1 - ( T + 1 )^{d p N_k}} \bigg) + 
\end{equation}
\[
+ \sum_{k = 1}^l \sum_{i = r + 1}^n \sum_{j = 0}^{( p-3 )/2} c_{i, j, k} \sum_{a = 1, p \not \mid a}^{dp} \bigg( {a^{2 j - 1} \chi_i ( a ) ( T + 1 )^{a N_k} \over 1 - ( T + 1 )^{d p N_k}} \bigg) \in \F_{q_1}[T].
\]
To finish the proof we shall prove that~(\ref{eqn:rel49}) is satisfied only if $c_{i, j, k} = 0$ for all $i, j, k$, obtaining a contradiction. 
 
We need the following remark:

\begin{rem}\label{rem:Osservazione11}
Let $ V $ be a vectorial space over a field $ F $ and $W$ a sub-space of $V$.
Moreover let $ \phi $ be an endomorphism of $ V $ such that $\phi ( W ) \subseteq W$.
We remark that if $ m $ is a positive integer and $v_1, v_2, \ldots , v_m \in V$ are eigen-vectors of $\phi$ with non-zero eigen-values $\lambda_1, \lambda_2, \ldots , \lambda_m$ such that $\lambda_i \neq \lambda_j$ if $i \neq j$ and if
\[
v_1 + v_2 + \ldots + v_m \in W,
\]   
then $v_i \in W$ for all $i$.
\end{rem}
  
If we apply Remark~\ref{rem:Osservazione11} in the particular case where $F = \F_{q_1}$, $V = \F_{q_1} ( T )$, $W = \F_{q_1} [T]$, $m = p - 1$, $\phi = \overline{D}$, $\lambda_b = b$ for all $1 \leq b \leq p - 1$ and
\[
v_b = V_b ( T ) = \sum_{k = 1}^l \sum_{i = 0}^r \sum_{j = 0}^{( p-3 )/2} c_{i, j, k} \sum_{a N_k \equiv b \mod ( p )} \bigg( {a^{2j} \chi_i ( a ) ( T + 1 )^{a N_k} \over 1 - ( T + 1 )^{d p N_k}} \bigg) + 
\]
\[
+ \sum_{k = 1}^l \sum_{i = r + 1}^n \sum_{j = 0}^{( p-3 )/2} c_{i, j, k} \sum_{a N_k \equiv b \mod ( p )}^{dp} \bigg( {a^{2 j - 1} \chi_i ( a ) ( T + 1 )^{a N_k} \over 1 - ( T + 1 )^{d p N_k}} \bigg),  
\]
by~(\ref{eqn:rel49}) we get $V_b ( T ) \in \F_{q_1} [T]$.    
Multiply $V_b ( T )$ by $1 - ( T + 1 )^{d p N_l}$.
Then
\[
( 1 - ( T + 1 )^{d p N_l} ) V_b ( T ) \in ( 1 - ( T + 1 )^{d p N_l} ) \F_{q_1}[T].
\]
Let us recall that $p$ does not divide $N_k$ for all $k$. 
Observe that if $\zeta$ is a primitive $d N_l$-th root of unity, then $\zeta - 1$ is a root of $( 1 - ( T + 1 )^{d p N_l} ) V_b ( T )$.
Since $N_l > N_k$ for all $k < l$, $\zeta - 1$ is a zero of
\[
( 1 - ( T + 1 )^{d p N_l} ) \sum_{k = 1}^{l-1} \sum_{i = 0}^r \sum_{j = 0}^{( p-3 )/2} c_{i, j, k} \sum_{a N_k \equiv b \mod ( p )} \bigg( {a^{2j} \chi_i ( a ) ( T + 1 )^{a N_k} \over 1 - ( T + 1 )^{d p N_k}} \bigg) + 
\]
\[
+ ( 1 - ( T + 1 )^{d p N_l} ) \sum_{k = 1}^{l-1} \sum_{i = r + 1}^n \sum_{j = 0}^{( p-3 )/2} c_{i, j, k} \sum_{a N_k \equiv b \mod ( p )}^{dp} \bigg( {a^{2 j - 1} \chi_i ( a ) ( T + 1 )^{a N_k} \over 1 - ( T + 1 )^{d p N_k}} \bigg).\]  
Then we get
\begin{equation}\label{eqn:relsola2}
\sum_{i = 1}^r \sum_{j = 0}^{( p-3 )/2} c_{i, j, l} \sum_{a N_l \equiv b \mod ( p )} a^{2 j} \chi_i ( a ) \zeta^{a N_l} + 
\end{equation}
\[
+ \sum_{i = r + 1}^n \sum_{j = 0}^{( p-3 )/2} c_{i, j, l} \sum_{a N_l \equiv b \mod ( p )} a^{2 j - 1} \chi_i ( a ) \zeta^{a N_l} = 0.    
\]
Observe that since $ p $ does not divide $d$, $\{ a, p + a, \ldots , p ( d - 1 ) + a \}$ is a set of representants of all the classes modulo $d$.
Moreover observe that, since $\zeta$ is a primitive $d N_l$-th root of unity, $\zeta^\prime = \zeta^{N_l}$ is a primitive $d$-th root of unity.
Let $k \in \Z/ p \Z$ such that $k N_l \equiv b \mod ( p )$.
We can rewrite~(\ref{eqn:relsola2}) as  
\begin{equation}\label{eqn:relsola3}
\sum_{i = 1}^r \sum_{j = 0}^{( p-3 )/2} c_{i, j, l} k^{2j} \sum_{h = 0}^{d - 1} \chi_i ( h ) \zeta^{\prime h} + \sum_{i = r + 1}^n \sum_{j = 0}^{( p-3 )/2} c_{i, j, l} k^{2j - 1} \sum_{h = 0}^{d - 1} \chi_i ( h ) \zeta^{\prime h} = 0.
\end{equation}
Then for all primitive $d$-th root of unity,~(\ref{eqn:relsola3}) must be satisfied. 

Set 
$$
x_{i, k} =
\begin{cases} 
\sum_{j = 0}^{( p - 3 )/2} c_{i, j, l} k^{2 j},&\;\hbox{if} \ 1 \leq i \leq r; \\
\sum_{j = 0}^{( p-3 )/2} c_{i, j, l} k^{2 j - 1}&\;\hbox{otherwise} 
\end{cases} 
$$
and let
\[
\{\zeta_1, \zeta_2, \ldots , \zeta_{\phi ( d )} \}
\]
be the set of primitive $d$-th roots of unity in $\overline{\F_p}$ (recall that $p$ does not divide $d$). 
Then by~(\ref{eqn:relsola3}) we have 
\begin{equation}\label{eqn:relsola4}
\sum_{i = 1}^n x_{i, k} \sum_{h = 0}^{d - 1} \chi_i ( h ) \zeta_c^h = 0.
\end{equation} 
for all $1 \leq c \leq \phi ( d )$. 
Hence we have a system of $n$ unknows ($x_{1, k} , \ldots , x_{n, k}$) and $\phi ( d )$ equations (one for all primitive $d$-th root of unity).
Observe that $n < \phi ( d )$.
Indeed by definition $n$ is less than the number of the characters of conductor $d \geq 2$ distinct $\mod ( \pi )$.
Since we have only $\phi ( d )$ characters whose conductor divides $d$ and since the trivial character has conductor $1 \neq d$, we have $n \leq \phi ( d ) - 1$.
Thus the number of equations of the system is greater then the number of its unknowns.  
We shall prove that the system has the unique solution $( 0, 0, \ldots , 0 )$.

Let $B$ be the matrix associated to the system~(\ref{eqn:relsola4}).
Then
$$
B = \left(
\begin{array}{cccc}
\sum_{h = 0}^{d - 1} \chi_1 ( h ) \zeta_1^h & \sum_{h = 0}^{d - 1} \chi_2 ( h ) \zeta_1^h & \cdots & \sum_{h = 0}^{d - 1} \chi_n ( h ) \zeta_1^h \\
\sum_{h = 0}^{d - 1} \chi_1 ( h ) \zeta_2^h & \sum_{h = 0}^{d - 1} \chi_2 ( h ) \zeta_2^h & \cdots & \sum_{h = 0}^{d - 1} \chi_n ( h ) \zeta_2^h \\
\vdots & \vdots & \ddots & \vdots \\
\sum_{h = 0}^{d - 1} \chi_1 ( h ) \zeta_{\phi ( d )}^h & \sum_{h = 0}^{d - 1} \chi_2 ( h ) \zeta_{\phi ( d )}^h & \cdots & \sum_{h = 0}^{d - 1} \chi_n ( h ) \zeta_{\phi ( d )}^h\\
\end{array}
\right).
$$
Observe that $B = C E$, where 
$$
C = \left(
\begin{array}{cccc}
1 & \zeta_1 & \cdots & \zeta_1^{d-1} \\
1 & \zeta_2 & \cdots & \zeta_2^{d-1} \\
\vdots & \vdots & \ddots & \vdots \\
1 & \zeta_{\phi ( d )} & \cdots & \zeta_{\phi ( d )}^{d - 1}\\
\end{array}
\right),
E = \left(
\begin{array}{cccc}
\chi_1 ( 0 ) & \chi_2 ( 0 ) & \cdots & \chi_n ( 0 ) \\
\chi_1 ( 1 ) & \chi_2 ( 1 ) & \cdots & \chi_n ( 1 ) \\
\vdots & \vdots & \ddots & \vdots \\
\chi_1 ( d - 1 ) & \chi_2 ( d - 1 ) & \cdots & \chi_n ( d - 1 ) \\
\end{array}
\right).
$$       
In the following lemma will prove that $\ker ( B )$ is equals to $\{ ( 0, 0, \ldots , 0 ) \}$, which implies $x_{1, k} = x_{2, k} = \ldots = x_{n, k} = 0$.

\begin{lem}\label{lem:lem11}
We have:
\begin{enumerate}
\item The rank of $C$ is $\phi ( d )$.\\
\item The set 
\[
\mathcal{B} = \{ ( 1, \zeta, \ldots, \zeta^{d-1} ), \ \zeta \in \mu_d, \ \zeta \ {\rm not} \ {\rm primitive} \}
\]
is a basis of $\ker ( C )$.\\
\item $\ker ( B ) = \{ ( 0, 0, \ldots , 0 ) \}$.  
\end{enumerate}
\end{lem}

\DIM 1. Let $C^\prime$ be the matrix whose columns coincide with the first $\phi ( d )$ columns of the matrix $C$.
Then $C^\prime$ is a square matrix equals to
$$
C^\prime = \left(
\begin{array}{cccc}
1 & \zeta_1 & \cdots & \zeta_1^{\phi ( d )-1} \\
1 & \zeta_2 & \cdots & \zeta_2^{\phi ( d ) -1} \\
\vdots & \vdots & \ddots & \vdots \\
1 & \zeta_{\phi ( d )} & \cdots & \zeta_{\phi ( d )}^{\phi ( d ) - 1}\\
\end{array}
\right).
$$ 
Observe that $C^\prime$ is a Vandermonde matrix. 
Thus its determinant is equals to 
\[
\prod_{1 \leq r < s \leq \phi ( d )} ( \zeta_r - \zeta_s ).
\]     
Since $p$ does not divide $d$, we have
\[
\zeta_r \neq \zeta_s. 
\] 
Thus the determinant of $C^\prime$ is non-zero.  
Therefore $C$ has rank $\phi ( d )$, since its square $\phi ( d ) \times \phi( d )$ sub-matrix $C^\prime$ has non-zero determinant. 

2. Since $C$ has rank $\phi ( d )$, the dimension of $\ker ( C )$ is $d - \phi ( d )$.
Observe that $\mathcal{B}$ has $d - \phi ( d )$ elements. 
Hence $\mathcal{B}$ is a basis of $\ker ( C )$ if and only if $\mathcal{B} \subseteq \ker ( C )$ and the elements of $\mathcal{B}$ are linearly independent.

First let us prove that $\mathcal{B} \subseteq \ker ( C )$.  
Let $\zeta \in \mu_d$ be a $d$-th root of unity which is not primitive.
Observe that $( 1, \zeta, \ldots , \zeta^{d-1} ) \in \ker ( C )$ if and only if
\[
\sum_{h = 0}^{d - 1} \zeta_i^h \zeta^h = 0
\]  
for all $i$.
Since
\[
\zeta_i \zeta \sum_{h = 0}^{d - 1} \zeta_i^h \zeta^h = \sum_{h = 0}^{d - 1} \zeta_i^h \zeta^h                
\]
and $\zeta_i \zeta \neq 1$ because $p$ does not divide $d$, it follows that
\[
\sum_{h = 0}^{d - 1} \zeta_i^h \zeta^h = 0.
\]
Thus $( 1, \zeta, \ldots , \zeta^{d-1} ) \in \ker ( C )$ for all non primitive $d$-th root of unity $\zeta$.
Then $\mathcal{B} \subseteq \ker ( C )$. 

Denote still by $\zeta$ an element of $\mu_d$ whose order is different from $d$. 
Set $\beta_\zeta \colon \Z / d \Z \rightarrow \mu_d$ the character which sends $i \in \Z /d \Z$ to $\zeta^i$.
Observe that if $\zeta^\prime \in \mu_d$ and $\zeta \neq \zeta^\prime$, 
\[
\beta_\zeta ( 1 ) = \zeta \neq \zeta^\prime = \beta_{\zeta^\prime} ( 1 ),
\] 
since $p$ does not divide $d$.
Thus the characters $\beta_\zeta$ are all distinct. 
Hence the theorem of the linear independence of characters imply that $\beta_\zeta$ are linearly independent over $\overline{\F_p}$.  
From this fact it follows that the vectors $( 1, \zeta, \ldots , \zeta^{d-1} )$ are linearly independent for all non primitive $d$-th root of unity $\zeta$.    
Hence $\mathcal{B}$ is a basis of $\ker ( C )$.

3. Using the previous notation for all non primitive $d$-th root of unity $\zeta$, let $\beta_\zeta$ be the function which sends $j \in \Z / d \Z$ to $\zeta^j$.         
Observe that every non trivial linear combination of the vectors $( \chi_i ( 0 ), \chi_i ( 1 ), \ldots , \chi_i ( d - 1 ) )$ for $1 \leq i \leq n$ is not in $\ker ( C )$ if and only if the functions $\chi_i$ and $\beta_\zeta$ for $1 \leq i \leq n$ and non primitive $d$-th root of unity $\zeta$ are linealry independent over $\overline{\F_p}$ (here $\chi_i$ is considered as a function of $\Z / d \Z$ over $\F_q \subseteq \overline{\F_p}$).
Suppose that such fonctions are dependent.
Then we can choose a minimal $r$, non-zero $\lambda_i$ and $\mu_\zeta$ in $\overline{\F_p}$ such that 
\begin{equation}\label{eqn:rel111}
\sum_{i = 1}^r \lambda_i \chi_i + \sum_{\zeta {\rm not} \ {\rm primitive}} \mu_\zeta \beta_\zeta = 0.
\end{equation}
First observe that $r \geq 2$.
Indeed if $r = 0$ then~(\ref{eqn:rel111}) would imply the linear dependence of the elements of $\mathcal{B}$ against 2.. 
Moreover if $r = 1$ then there would exist a character $\chi_i$ of conductor $d$ such that $( \chi_i ( 0 ), \chi_i ( 1 ), \ldots , \chi_i ( d - 1 ) ) \in \ker ( C )$.  
Then we would have
\[
\sum_{h = 0}^{d-1} \zeta_j^h \chi_i ( h ) = 0
\]
for all primitive $d$-th root of unity $\zeta_j$, which contradicts~\cite[Lemma 4.8]{4}.
             
Let $b \in ( \Z / d \Z )^\ast$ such that $\chi_1 ( b ) \neq \chi_2 ( b )$ (such $b$ exists because recall that, by hypothesis, $\chi_1$ is different from $\chi_2$ modulo $( \pi )$).   
Hence by~(\ref{eqn:rel111}), for all $z \in \Z / d \Z$ we have
\begin{equation}\label{eqn:rel112}
\sum_{i = 1}^r \lambda_i \chi_i ( b z ) + \sum_{\zeta {\rm not} \ {\rm primitive}} \mu_\zeta \beta_\zeta ( b z ) = 0.
\end{equation}
Observe that the function which sends $z \in \Z / d \Z$ to $\beta_\zeta ( b z )$ coincides with the function $\beta_{\zeta^b}$. 
Hence we can rewrite~(\ref{eqn:rel112}) as 
\begin{equation}\label{eqn:rel113}
\sum_{i = 1}^r \lambda_i \chi_1 ( b ) \chi_i ( z ) + \sum_{\zeta {\rm not} \ {\rm primitive}} \chi_1 ( b ) \mu_\zeta \beta_{\zeta^b} ( z ) = 0.
\end{equation}
If we multiply~(\ref{eqn:rel111}) by $\chi_1 ( b )$ and we subtract it to~(\ref{eqn:rel113}),we get a non trivial relation with less than $r$ characters $\chi_i$ and it is impossible by the minimality of $r$.

Finally consider the matrix $B$.
Let $v = ( \lambda_1, \lambda_2, \ldots , \lambda_n ) \in \ker ( B )$.
Since $B = CE$, then $E ( v ) \in \ker ( C )$.
This fact implies that 
\[
\sum_{ i = 1 }^n \lambda_i ( \chi_i ( 0 ), \chi_i ( 1 ), \ldots , \chi_i ( d - 1 ) ) \in \ker( C ).
\]
But we have previously proved that this relation is possible only if $\lambda_i = 0$ for all $i$.
Thus $v = ( 0, 0, \ldots , 0 )$. 
\CVD

By the previous lemma we immediately get $x_{i, k} = 0$ for all $1 \leq i \leq n$ and for all $1 \leq k \leq p-1$. 
Remember that, by definition,
$$
x_{i, k} =
\begin{cases} 
\sum_{j = 0}^{( p - 3 )/2} c_{i, j, l} k^{2 j},&\;\hbox{if} \ 1 \leq i \leq r; \\
\sum_{j = 0}^{( p-3 )/2} c_{i, j, l} k^{2 j - 1}&\;\hbox{otherwise.} 
\end{cases} 
$$ 
We shall prove that the relation $x_{i, k} = 0$ for all $i, k$ implies $c_{i, j, l} = 0$ for all $i, j$.
We just consider the case $i \leq r$ (the proof in the other case is very similar). 
Let $i$ be an integer between $1$ and $r$ and set $c_{i, j, l} = y_j$. 
Since $x_{i, k} = 0$ for all $k$ between $1$ and $p-1$, we have the following relations:  
\begin{equation}\label{eqn:rel214}
\sum_{j = 0}^{( p-3 )/2} y_j k^{2 j} = 0.
\end{equation} 
The matrix $M$ associated to the first $( p - 1 )/2$ equations of the system~(\ref{eqn:rel214}) is given by:
$$
M =
\left(
\begin{array}{cccc}
1 & 1 & \cdots & 1 \\
1 & 4 & \cdots & 4^{( p - 3 )/ 2} \\
\vdots & \vdots & \ddots & \vdots \\
1 & {( p-1 )^2 \over 4} & \cdots & \bigg( {( p-1 )^2 \over 4} \bigg)^{( p - 3 )/2}\\
\end{array}
\right) =
\left(
\begin{array}{cccc}
1 & \alpha_1 & \cdots & \alpha_1^{( p-3 )/2} \\
1 & \alpha_2 & \cdots & \alpha_2^{( p-3 )/2} \\
\vdots & \vdots & \ddots & \vdots \\
1 & \alpha_{( p - 1 )/2} & \cdots & \alpha_{( p-1 )/2}^{( p - 3 ) /2}\\
\end{array}
\right)                  
$$     
So $M$ is a Vandermonde matrix and its determinant is equals to:
\[
{\rm det} ( M ) = \prod_{1 \leq r < s \leq ( p-1 )/2} ( \alpha_r - \alpha_s ) \neq 0.
\]  
It follows that the only solution of the system~(\ref{eqn:rel214}) is $y_j = 0$ for all $j$.

Then we have proved that the coefficients $c_{i, j, l} = 0$ for all $i$ and $j$.
Since $N_{l-1} > N_k$ for all $k < l-1$, if we replace $l$ with $l-1$ with the same procedure we can prove that the coefficients $c_{i, j, l-1} = 0$ for all $i, j$ and so on.
Thus $c_{i, j, k} = 0$ for all $i, j, k$, which is a contradiction.
\CVD

\section{The case $d = 1$}\label{sec:partequarta}

The aim of this section is to prove Theorem~\ref{teo:teo1} in the case where
$d = 1$.
In other words, let $\chi$ the trivial character. 
We shall prove that 
\[
{\rm dim}_\Omega ( \Omega + \Omega \overline{f ( T, \chi \omega^0 )} + \Omega \overline{f ( T, \chi \omega^2 )} + \ldots + \Omega \overline{f ( T, \chi \omega^{p-3} ) )} = {p + 1 \over 2}.   
\]
As we have already remarked is Section~\ref{sec:parteprima} we shall modify the proof of the case $d \geq 2$. 
Let us give a preliminar reason for this.
Let $\delta \in \Z/ ( p-1 ) \Z$ be odd.
Then by Lemma~\ref{lem:lem4}(6.) we have
\[
f ( ( T + 1 )^{-1} - 1, \chi \omega^{\delta + 1} ) = \Gamma_\delta \gamma_{-\delta} U ( ( F_\chi ( T ) ).
\]
Since $\chi$ is the trivial character we have
\begin{equation}\label{eqn:relmontagna1}
F_\chi ( T ) = \sum_{a = 0}^{p - 1} { ( T + 1 )^a \over 1 - ( T + 1 )^p} - 1 = - {1 \over T} - 1.
\end{equation}
Thus $\overline{F_\chi} ( T ) \not \in \F_p[[T]]$ and we shall see that this fact does not allow us to apply Proposition~\ref{pro:pro2} (observe that if $\chi^\prime$ is not the trivial character then $\overline{F_{\chi^\prime}} ( T ) \in \F_q[[T]]$ and  see also Remark~\ref{rem:Osservazione2}).
The following lemma explains how we can solve this problem.

\begin{lem}\label{lem:lem6}
Let 
\[
\widetilde{F_\chi} ( T ) = ( p + 1 ) F_\chi ( ( T + 1 )^{p+1} - 1 ) - F_\chi ( T ).
\]  
Then $\widetilde{F_\chi} ( T ) \in \Z_p[[T]]$.

Moreover 
\[
T \overline{f ( ( T + 1 )^{-1} - 1, \chi \omega^{\delta + 1} )} = \overline{\Gamma}_\delta \overline{\gamma}_{-\delta} \overline{U} ( \overline{\widetilde{F_\chi}} ( T ) )
\]
for all odd $\delta \in \Z/ ( p-1 ) \Z$.

Finally $1, T \overline{f ( ( T + 1 )^{-1} - 1, \chi \omega^0 )}, \ldots , T \overline{f ( ( T + 1 )^{-1} - 1, \chi \omega^{p-3} )}$ are independent over $\Omega$ if and only if $1, \overline{f ( T , \chi \omega^0 )}, \ldots , \overline{f ( T, \chi \omega^{p-3} )}$ are independent over $\Omega$.    
\end{lem}

\DIM Observe that
\begin{align*}
\widetilde{F_\chi} ( T ) & = ( p + 1 ) F_\chi ( ( T + 1 )^{p+1} - 1 ) - F_\chi ( T ) \\ 
& = -{( p + 1 )  \over ( T + 1 )^{p + 1} - 1} + {1 \over T} - p \\
& = { ( T + 1 )^{p+1} - ( p + 1 ) T - 1 \over T ( ( T + 1 )^{p+1} - 1 ) } - p.
\end{align*}
We remark that $ T^2 $ exactly divides $T ( ( T + 1 )^{p+1} - 1 )$.
Moreover $T^2$ divides $( T + 1 )^{p + 1} - ( p + 1 ) T - 1$.
Then we immediately get $\widetilde{F_\chi} ( T ) \in \Z_p[[T]]$.     

Let $L$ be a finite extension of $\Q_p$ and let $\O_L$ be its valuation ring.
Let $d$ be an integer not divided by $p$ and set $\kappa_0 = 1 + dp$. 
Remember that in Section~\ref{sec:parteseconda} we have defined an action of $\O_L[[t]]$ over $\O_L[[T]]$ such that, for all $a \in \Z_p$ and $F ( T ) \in \O_L[[T]]$,
\[
( t + 1 )^a F ( T ) = F ( ( T + 1 )^{\kappa_0^a} - 1 ).
\] 
Set $L = \Q_p$, $d = 1$ and $\kappa_0 = 1 + p$.
Then, since
\[
\widetilde{F}_\chi ( T ) = F_\chi ( ( T + 1 )^{p+1} - 1 ) - F_\chi ( T ) + p F_\chi ( ( T + 1 )^{p+1} - 1 ), 
\]
we have
\[
\widetilde{F_\chi} ( T ) \equiv t F_\chi ( T ) \mod ( p ).
\]      
By Lemma~\ref{lem:lem5}(2.) we get
\[
\overline{\Gamma_\delta} \overline{\gamma_{-\delta}} \overline{U} ( ( \overline{\widetilde{F_\chi}} ( T ) ) = \overline{\Gamma_\delta} \overline{\gamma_{-\delta}} \overline{U} ( t \overline{F_\chi} ( T ) ) = T \overline{\Gamma_\delta} \overline{\gamma_{-\delta}} \overline{U} ( \overline{F_\chi} ( T ) ). 
\]  
Then, by Lemma~\ref{lem:lem4}(6.), we immediately get
\[
T \overline{f ( ( T + 1 )^{-1} - 1, \chi \omega^{\delta + 1} )} = \overline{\Gamma}_\delta \overline{\gamma}_{-\delta} \overline{U} ( \overline{\widetilde{F_\chi}} ( T ) ).      
\]

Finally $1, T \overline{f ( ( T + 1 )^{-1} - 1, \chi \omega^0 )}, \ldots , T \overline{f ( ( T + 1 )^{-1} - 1, \chi \omega^{p-3} )}$ are independent over $\Omega$ if and only if  $1, \overline{f ( T, \chi \omega^0 )}, \ldots , \overline{f ( T, \chi \omega^{p-3} )}$ are independent over $\Omega$, since $T \in \Omega$ and the image of $\Omega$ via the endomorphism of $\overline{\F_p} ( ( T ) )$ which sends $F ( T ) \in \overline{\F_p} ( ( T ) )$ in $F ( ( T + 1 )^{-1} - 1 )$ is $\Omega$. 
\CVD  

By the previous lemma to finish the proof of Theorem~\ref{teo:teo1} it suffices to show that $1, T \overline{f ( ( T + 1 )^{-1} - 1, \chi \omega^0 )}, \ldots , T \overline{f ( ( T + 1 )^{-1} - 1, \chi \omega^{p-3} )}$ are linearly independent over $\Omega$. 
Since, always by the previous lemma, for all odd $\delta \in \Z / ( p-1 ) \Z$ we have
\[
T \overline{f ( ( T + 1 )^{-1} - 1, \chi \omega^{\delta + 1} )} = \overline{\Gamma}_\delta \overline{\gamma}_{-\delta} \overline{U} ( \overline{\widetilde{F_\chi}} ( T ) )  
\]
and $\overline{\widetilde{F_\chi}} ( T ) \in \F_p[[T]]$, to show the linear independence we can easily adapt our proof of Theorem~\ref{teo:teo1} in the case $d \geq 2$. 
\vspace{5 pt}

{\bf Proof of Theorem~\ref{teo:teo1} in the case $d = 1$.}            
First observe that by Lemma~\ref{lem:lem6} to prove the assumption it suffices to show that 
\[
1, T \overline{f ( ( T + 1 )^{-1} - 1 , \chi \omega^0 )}, \ldots , T \overline{f ( ( T + 1 )^{-1} - 1, \chi \omega^{p-3} )}
\]
are linearly independent over $\Omega$.
Suppose that this is not the case.
Then there exist $h_0 ( T ), \ldots , h_{( p-3 ) / 2} ( T ) \in \Omega$ such that
\begin{equation}\label{eqn:rel31}
\sum_{j = 0}^{( p-3 )/2} h_j ( T ) T \overline{f\bigg( {1 \over T + 1} - 1, \omega^{2j} \bigg)} \in \Omega,
\end{equation}     
with $h_j ( T ) \neq 0$ for a certain $j$.
Observe that without loss of generality we can suppose that $h_{j} ( T ) \in A$ for all $i$ and that   
\[
\sum_{j = 0}^{( p-3 )/2} h_j ( T ) T \overline{ f\bigg( {1 \over T + 1} - 1, \omega^{2j} \bigg)} \in A.
\]
By Lemma~\ref{lem:lem6} we get 
\[
\sum_{j = 0}^{( p-3 )/2} h_j ( T ) \overline{\Gamma}_{2j - 1} \overline{\gamma}_{- 2j + 1 } \overline{U} ( \overline{\widetilde{F_\chi}} ( T ) ) \in A.
\]

By Lemma~\ref{lem:lem3}(3.) and since $\overline{D} \overline{\gamma_{\delta}} = \overline{\gamma_{\delta + 1}} \overline{D}$ for all $\delta \in \Z/ ( p-1 ) \Z$, we have
\[
\sum_{j = 0}^{( p-3 )/2} h_j ( T ) \overline{\Gamma}_{- 1} \overline{\gamma}_1 \overline{D}^{2j} \overline{U} ( \overline{\widetilde{F_\chi}} ( T ) ) \in A.
\]     
Moreover, applying Lemma~\ref{lem:lem5}, we get
\[
\sum_{j = 0}^{( p-3 )/2} \overline{\Gamma}_{- 1} \overline{\gamma}_1 \overline{D}^{2j} \overline{U} ( h_j ( t ) \overline{\widetilde{F_\chi}} ( T ) ) \in A.
\]
From Lemma~\ref{lem:lem3}(2.) it follows 
\begin{equation}\label{eqn:relmontagna2} 
\sum_{j = 0}^{( p-3 )/2}  \overline{\gamma}_1 \overline{D}^{2j} \overline{U} ( h_j ( t ) \overline{\widetilde{F_\chi}} ( T ) ) \in A.
\end{equation}

For all $j$ such that $0 \leq j \leq ( p - 3 )/2$ set
\[
F_j ( T ) = \overline{D}^{2j} \overline{U} ( \overline{\widetilde{F_\chi}} ( T ) ). 
\] 
Then we can rewrite~(\ref{eqn:relmontagna2}) in the following way:
\begin{equation}\label{eqn:relmontagna3} 
\sum_{j = 0}^{( p - 3 )/2} \overline{\gamma}_1 ( h_j ( t ) F_j ( T ) ) \in A.
\end{equation}

By Lemma~\ref{lem:lem4}(4.) we have
\[
F_\chi ( ( T + 1 )^{-1} - 1 ) = -F_\chi ( T ) - 1.
\]
Then
\begin{align*}
\overline{\widetilde{F_\chi}} ( ( T + 1 )^{-1} - 1 ) & = \overline{F_\chi} ( ( T + 1 )^{- ( p + 1 )} - 1 ) - \overline{F_\chi} ( ( T + 1 )^{-1} - 1 ) \\
& = - \overline{F_\chi} ( ( T + 1 )^{p + 1} - 1 ) - 1 + \overline{F_\chi} ( T ) + 1 \\
& = - \overline{\widetilde{F_\chi}} ( T ).
\end{align*}   
Moreover observe that for all $j$ between $0$ and $( p - 3 )/2$, 
\[
( \overline{D}^{2 j} \overline{U} ( \overline{\widetilde{F_\chi}} ) ) ( ( T + 1 )^{-1} -1 ) = \overline{D}^{2j} \overline{U} ( \overline{\widetilde{F_\chi}} ( ( T + 1 )^{-1} - 1 ).
\] 
Then for all $j$
\[ 
( \overline{D}^{2 j} \overline{U} ( \overline{\widetilde{F_\chi}} ) ) ( ( T + 1 )^{-1} -1 ) = - \overline{D}^{2j} \overline{U} ( \overline{\widetilde{F_\chi}} ( T ) ).
\]
It follows that
\begin{equation}\label{eqn:relmontagna4} 
F_j ( T ) = - F_j ( ( T + 1 )^{-1} - 1 )
\end{equation}  
for all $j$.

Consider the least field which contains $\F_p$ and all the coefficients of $h_ j ( T )$ for all $j$.
Since $h_j ( T ) \in A \subseteq \overline{\F_p}[[T]]$, such field is a finite extension of $\F_p$.
Call it $\F_{q_1}$ and write $h_j ( t ) = \sum_{b \in \Z_p} c_{j, b} ( t + 1 )^b$ with $c_{j, b} \in \F_{q_1}$.
Moreover observe that since $h_j ( T ) \neq 0$ for certain integer $j$, there exist $j, b$ such that $c_{j, b} \neq 0$.
Set 
\[
G_b ( T ) = \sum_{j = 0}^{( p - 3 )/2} c_{j, b} F_j ( T )
\]
and observe that $G_b ( T ) \in \F_{q_1}[[T]] \cap \F_{q_1} ( T )$ for all $b$. 
Since 
\[
( t + 1 )^b ( G_b ( T ) ) = G_b ( ( T + 1 )^{\kappa_0^b} - 1 ), 
\]
by~(\ref{eqn:relmontagna3}) we have
\begin{equation}\label{eqn:relmontagna6}
\overline{\gamma}_1 ( \sum_{b \in \Z_p}  G_b ( ( T + 1 )^{\kappa_0^b} - 1 ) ) \in A.
\end{equation}
   
Choose a subset of $\mu_{p-1}$ whose elements represent all the classes of the group $\mu_{p - 1} / \{-1, 1\}$ and call it $S$.
From~(\ref{eqn:relmontagna4}) it follows that 
\begin{equation}\label{eqn:relmontagna7}
G_b ( T ) =  - G_b ( ( T + 1 )^{-1} - 1 )
\end{equation} 
for all $b$.
Thus by~(\ref{eqn:relmontagna6}) we get
\[
\sum_{\eta \in S} \sum_{b \in \Z_p} \eta G_b ( ( T + 1 )^{\eta \kappa_0^b} - 1 ) \in A.
\]
Since $G_b ( T ) = 0$ for all but finitely many $b \in \Z_p$, there exists an integer $u$ such that 
\[
\sum_{\eta \in S} \sum_{k = 1}^u \eta G_{b_k} ( ( T + 1 )^{\eta \kappa_0^{b_k}} - 1 ) \in A.
\]
Moreover set $G_k ( T ) = G_{b_k} ( T )$.
Since $G_k ( T ) \in \F_{q_1}[[T]] \cap \F_{q_1} ( T )$ for all $k$, we can apply Proposition~\ref{pro:pro2}, obtaining that there exist an integer $l \leq u$, $b_1, b_2, \ldots b_l \in \Z_p$, $b_i \neq b_j$ for $i \neq j$, $\eta_1, \eta_2, \ldots \eta_l \in \mu_{p - 1}$ with $\eta_i \kappa_0^{b_i} \sim_{\Q^\ast} \eta_j \kappa_0^{b_j}$ for all $i$, $j$ and $\eta_i \kappa_0^{b_i} \neq \eta_j \kappa_0^{b_j}$ for $i \neq j$ such that   
\begin{equation}\label{eqn:relmontagna8}
\sum_{k = 1}^l \eta_k  G_k ( ( T + 1 )^{\eta_k \kappa_0^{b_k}} - 1 ) \in A.   
\end{equation}

\begin{rem}\label{rem:Osservazione2}
We remark that that the fact that for all $k$, $G_k ( T ) \in \F_{q_1}[[T]] \cap \F_{q_1} ( T )$ is necessary to apply Proposition~\ref{pro:pro2}.
This relation is verified since for all $j$, we have $F_j ( T ) \in \F_p [[T]] \cap \F_p ( T )$ which is an immediate consequence of the fact that $\overline{\widetilde{F_\chi}} ( T ) \in \F_p[[T]] \cap \F_p ( T )$.

Finally observe that in the case $d \geq 2$, for all character $\chi^\prime$ of conductor $d$ we have $\overline{F_{\chi^\prime}} ( T ) \in \F_q[[T]] \cap \F_q ( T )$. 
For this reason it is not necessary <<to perturbate>> $F_{\chi^\prime} ( T )$ to apply Proposition~\ref{pro:pro2}.
\end{rem}  

For all $1 \leq k \leq l$ write
\[
\eta_k \kappa_0^{b_k} = \eta_1 \kappa_0^{b_1} x_k,
\]  
where $x_k \in \Q^\ast \cap \Z_p^\ast$ and $x_i \neq x_k$ if $i \neq k$.
Recall that by~(\ref{eqn:relmontagna7}), we have $G_k ( T ) = - G_k ( ( T + 1 )^{-1} - 1 )$.
Hence we can suppose that $x_k > 0$ for all $k$. 
By~(\ref{eqn:relmontagna6}) we get
\[
\sum_{k = 1}^l \eta_k  G_k ( ( T + 1 )^{x_k} - 1 ) \in A.
\] 
Therefore there exist some integers $N_1 , N_2 , \ldots N_l$ not divided by $p$, such that $1 \leq N_1 < N_2 < \ldots < N_l$ and 
\[
\sum_{k = 1}^l  \eta_k G_k ( ( T + 1 )^{N_k} - 1 ) \in A.    
\]
By definition of $G_k ( T )$ this last relation becomes:
\begin{equation}\label{eqn:relmontagna9}
\sum_{k = 1}^l \eta_k \sum_{j = 0}^{( p-3 )/2} c_{j, k} F_j ( ( T + 1 )^{N_k} - 1 ) \in A.  
\end{equation}   

Now we want to compute $F_j ( T )$ for all $j$.     
First remember that
\[
F_\chi ( T ) = \sum_{a = 0}^{p-1} {( T + 1 )^a \over 1 - ( T + 1 )^p} - 1.
\] 
Then
\[
\overline{D}^{2j} ( \overline{F_\chi} ( T ) ) = \sum_{a = 1}^{p-1} {a^{2 i} ( T + 1 )^a \over 1 - ( T + 1 )^p}
\]
for all $j$ between $0$ and $( p-1 )/2$.
Since $\overline{U} = \overline{D}^{p-1}$ we have  
\[
\overline{D}^{2i} \overline{U} ( \overline{F_\chi} ( T ) ) = \sum_{a = 1}^{p-1} {a^{2 i} ( T + 1 )^a \over 1 - ( T + 1 )^p}.
\]
Remember that 
\[
\overline{\widetilde{F_\chi}} ( T ) = \overline{F_\chi} ( ( T + 1 )^{p+1} - 1 ) - \overline{F_\chi} ( T ) = t \overline{F_\chi} ( T ).
\]  
Then 
\begin{align*}
\overline{D}^{2j} \overline{U} ( \overline{\widetilde{F_\chi}} ( T ) ) & = t \overline{D}^{2j} \overline{U} ( \overline{F_\chi} ( T ) ) \\
& = t \bigg( \sum_{a = 1}^{p-1} {a^{2 i} ( T + 1 )^a \over 1 - ( T + 1 )^p} \bigg) \\
& = \sum_{a = 1}^{p - 1} {a^{2i} ( T + 1 )^{a ( p + 1 )} \over 1 - ( T + 1 )^{p ( p + 1 )}} - \sum_{a = 1}^{p - 1} {a^{2i} ( T + 1 )^a \over 1 - ( T + 1 )^p}.   
\end{align*}
Replacing in~(\ref{eqn:relmontagna9}) we get   
\[
\sum_{k = 1}^l \eta_k \sum_{j = 0}^{( p-3 )/2} c_{j, k} \sum_{a = 1}^{p - 1} \bigg( {a^{2j} ( T + 1 )^{a N_k ( p + 1 )} \over 1 - ( T + 1 )^{p N_k ( p + 1 )}} - {a^{2j} ( T + 1 )^{a N_k} \over 1 - ( T + 1 )^{p N_k}} \bigg) \in A.
\] 
Since by~\cite[Lemma 3.5]{1} a rational function $H ( T ) \in A$ if and only if there exists an integer $n$ such that $( T + 1 )^n H ( T )$ is a polynomial, we get
\begin{equation}\label{eqn:relmontagna10}
\sum_{k = 1}^l \eta_k \sum_{j = 0}^{( p-3 )/2} c_{j, k} \sum_{a = 1}^{p - 1} \bigg( {a^{2j} ( T + 1 )^{a N_k ( p + 1 )} \over 1 - ( T + 1 )^{p N_k ( p + 1 )}} - {a^{2j} ( T + 1 )^{a N_k} \over 1 - ( T + 1 )^{p N_k}} \bigg) \in \F_{q_1} [T]
\end{equation}
(recall that, by definition, $\F_{q_1}$ is the extension of $\F_p$ generated by $c_{j, k}$ for all $j, k$).

We apply Remark~\ref{rem:Osservazione11} in the particular case where $F = F_{q_1}$, $V = \F_{q_1} ( T )$, $W = \F_{q_1} [T]$, $m = p - 1$, $\phi = \overline{D}$, $\lambda_b = b$ for all $1 \leq b \leq p - 1$ and
\[
v_b = V_b ( T ) : = \sum_{k = 1}^l \eta_k \sum_{j = 0}^{( p-3 )/2} c_{j , k} \bigg( {a_{b, k}^{2j} ( T + 1 )^{a_{b,k} N_k ( p + 1 )} \over 1 - ( T + 1 )^{p N_k ( p + 1 )}} - {a_{b, k}^{2j} ( T + 1 )^{a_{b,k} N_k} \over 1 - ( T + 1 )^{p N_k}} \bigg),  
\]
where $a_{b, k}$ satisfies the relation $a_{b, k} N_k \equiv b \mod ( p )$.
Then by Remark~\ref{rem:Osservazione11}, $V_b ( T ) \in \F_{q_1}[T]$. 
Let us recall that $p$ does not divide $N_k$ for all $k$ between $1$ and $l$.
Since $p$ does not divide $b$, we have  $a_{b, k} \in \F_p^\ast$ for all $b, k$. 
 
Let $\zeta$ be a primitive $( p + 1 ) N_l$th root of unity and multiply $V_b ( T )$ by the polynomial $1 - ( T + 1 )^{p ( p + 1 ) N_l}$. 
We get
\begin{equation}\label{eqn:rel311}
Q( T ) + \eta_l \sum_{j = 0}^{( p - 3 )/2} c_{j, l} a_{b,l}^{2 j} ( T + 1 )^{p ( p + 1 ) N_l} = P( T ),
\end{equation} 
where $Q ( T ) \in \F_{q_1} ( T ), P ( T ) \in \F_{q_1}[T]$, $Q ( \zeta - 1 ) = 0$ (since $N_l > N_k$ for all $k < l$) and $P ( \zeta - 1 ) = 0$. 
Then if~(\ref{eqn:rel311}) is satisfied we have
\[
\eta_l \zeta^{a_b N_l ( p + 1 )} \sum_{j = 0}^{( p - 3 )/2} c_{j, l} a_{b, l}^{2 j} = 0 
\]
for all $b$ between $1$ and $p - 1$.
Since $\eta_l \zeta^{a_b N_l ( p + 1 )} = \eta_l \in \F_p^\ast$, we have
\begin{equation}\label{eqn:relmontagna11}
\sum_{j = 0}^{( p - 3 )/2} c_{j, l} a_{b,l}^{2 j} = 0, 
\end{equation}  
for all $b$. 
Observe that, since $1 \leq b \leq p-1$, $a_{b, l} N_l \equiv b \mod ( p )$ and $p$ does not divide $N_l$, for all $c \in \F_p^\ast$ there exists $b$ such that $a_{b,l} = c$.    
Then to find the $c_{j, l}$ satisfying~(\ref{eqn:relmontagna11}) is equivalent to find all the solutions of the system
\begin{equation}\label{eqn:rel312}
\sum_{j = 0}^{ ( p - 3 )/2} x_j k^{2 j} = 0,
\end{equation}         
where $x_j$ are the unknows of the system and $1 \leq k \leq p-1$. 
Observe that the system~(\ref{eqn:rel312}) is identical to the system~(\ref{eqn:rel214}) that we have studied at the end of Section~\ref{sec:parteterza}. 
Since we have already remarked that~(\ref{eqn:rel214}) has only the solution $( 0, \ldots , 0 )$, also~(\ref{eqn:rel312}) has the unique solution $( 0, \ldots , 0 )$.
This fact implies $c_{j, l} = 0$ for all $j$.

Since $N_{l-1} > N_k$ for all $k < l-1$, if we replace $l$ with $l-1$ with the same procedure we can prove that the coefficients $c_{j, l-1} = 0$ for all $j$ and so on.
Thus $c_{j, k} = 0$ for all $j, k$, which implies $h_j ( T ) = 0$ for all $j$.
Since we have supposed that there exists $j$ such that $h_j ( T ) \neq 0$, we obtain a contraddiction.     
\CVD


\begin{thebibliography}{[Bo--Gi--SO]}

\bibitem[1]{1} B. Anglès, {\it On the $p$-adic Leopoldt transform of a power series}, Acta Arith. {\bf 134}, {349-367}, (2008).

\bibitem[2]{2} S. Lang, {\it Cyclotomic field I, II}, Springer, 1990.

\bibitem[3]{3} W. Sinnott, {\it On the power series attached to $p$-adic $L$-functions}, J. Reine Angew. Math. {\bf 382}, {22-34}, (1987).

\bibitem[4]{4} L. C. Washinghton, {\it Introduction to cyclotomic fields}, 2nd ed., Springer, 1997.












\end{thebibliography}
\end{document}